\catcode`\@=11
\magnification 1200
\hsize=140mm \vsize=200mm
\hoffset=-4mm \voffset=-1mm
\pretolerance=500 \tolerance=1000 \brokenpenalty=5000

\catcode`\;=\active
\def;{\relax\ifhmode\ifdim\lastskip>\z@
\unskip\fi\kern.2em\fi\string;}

\catcode`\:=\active
\def:{\relax\ifhmode\ifdim\lastskip>\z@\unskip\fi
\penalty\@M\ \fi\string:}

\catcode`\!=\active
\def!{\relax\ifhmode\ifdim\lastskip>\z@
\unskip\fi\kern.2em\fi\string!}\catcode`\?=\active
\def?{\relax\ifhmode\ifdim\lastskip>\z@
\unskip\fi\kern.2em\fi\string?}

\frenchspacing

\newif\ifpagetitre        \pagetitretrue
\newtoks\hautpagetitre    \hautpagetitre={\hfil}
\newtoks\baspagetitre     \baspagetitre={\hfil}
\newtoks\auteurcourant    \auteurcourant={\hfil}
\newtoks\titrecourant     \titrecourant={\hfil}
\newtoks\hautpagegauche   \newtoks\hautpagedroite
\hautpagegauche={\hfil\tensl\the\auteurcourant\hfil}
\hautpagedroite={\hfil\tensl\the\titrecourant\hfil}

\newtoks\baspagegauche
\baspagegauche={\hfil\tenrm\folio\hfil}
\newtoks\baspagedroite
\baspagedroite={\hfil\tenrm\folio\hfil}

\headline={\ifpagetitre\the\hautpagetitre
\else\ifodd\pageno\the\hautpagedroite
\else\the\hautpagegauche\fi\fi}

\footline={\ifpagetitre\the\baspagetitre
\global\pagetitrefalse
\else\ifodd\pageno\the\baspagedroite
\else\the\baspagegauche\fi\fi}

\font\twbf=cmbx12\font\sc=cmcsc10

\font\tenhelv=phvr at 10pt
\font\sevenhelv=phvr at 7pt
\font\fivehelv=phvr at 5pt
\font\ethelv=phvr at 15pt
\font\tengoth=eufm10
\font\sevengoth=eufm7
\font\fivegoth=eufm5
\font\tenbb=msbm10
\font\sevenbb=msbm7
\font\fivebb=msbm5
\newfam\gothfam
\newfam\bbfam
\textfont\gothfam=\tengoth
\scriptfont\gothfam=\sevengoth
\scriptscriptfont\gothfam=\fivegoth

\textfont\bbfam=\tenbb
\scriptfont\bbfam=\sevenbb
\scriptscriptfont\bbfam=\fivebb

\def\date{{\the\day}\
\ifcase\month\or Janvier\or \F\'evrier\or Mars\or Avril
\or Mai\or Juin\or Juillet\or Ao\^ut\or Septembre
\or Octobre\or Novembre\or D\'ecembre\fi\ {\the\year}}

\def\cf{{\it cf.\/}\ }    \def\ie{{\it i.e.\/}\ }
 \def\up#1{\raise 1ex\hbox{\sevenrm#1}}
\def\cqfd{\unskip\kern 6pt\penalty 500
\raise -2pt\hbox{\vrule\vbox to 10pt{\hrule width 4pt\vfill
\hrule}\vrule}\par}
\catcode`\@=12

\def\sn{\nobreak\smallskip}

\def\PI{{\mit\Pi}} 

\def\ref#1&#2&#3&#4&#5\par{\par{\leftskip = 5em{\noindent
\kern-5em\vbox{\hrule height0pt depth0pt width
5em\hbox{\bf[\kern2pt#1\unskip\kern2pt]\enspace}}\kern0pt}
{\sc\ignorespaces#2\unskip},\
{\rm\ignorespaces#3\unskip}\
{\sl\ignorespaces#4\unskip\/}\
{\rm\ignorespaces#5\unskip}\par}}

\def\exo#1{\goodbreak\vskip 12pt plus 20pt minus 2pt%
\line{\noindent\hss\bf%
\uppercase\expandafter{\romannumeral#1}\hss}\nobreak\vskip
12pt }
\def \titre#1\par{\null\vskip
1cm\line{\hss\vbox{\let\\=\cr\twbf\halign%
{\hfil##\hfil\crcr#1\crcr}}\hss}\vskip 1cm}

\def\frac#1#2{{#1\over#2}}

\def\comp{\;{}^{{}_\vert}\!\!\!{\rm C}}

\def\nat{{{\rm I}\!{\rm N}}}

\def\rat{{\rm Q}\kern-.65em {}^{-{}_/}}

\def\N#1{\left\Vert#1\right\Vert}

\let\isspt=\i
\def\i{\ifmmode\infty\else\isspt\fi}

\def\dess#1by#2(#3){\vbox to #2{\hrule width #1
height 0pt depth 0pt\vfill\special{picture #3}}}
\def\Dess#1by#2(#3 scaled#4){{\dimen10=#1\dimen11=#2
\divide\dimen10 by1000\multiply\dimen10 #4
\divide\dimen11 by1000\multiply\dimen11 #4
\vbox to\dimen11{\hrule width\dimen10
height\z@ depth\z@\vfill\special{picture #3 scaled#4}}}}

\def\equipe{\setbox10=\hbox{\dess36mm by 26mm(analyse
scaled 800)}
\setbox12 = \vbox {\hsize =8cm
\centerline{\ethelv EQUIPE D'ANALYSE}
\vskip-4pt
\centerline{\fivehelv URA 754 -- CNRS}
\centerline{\tenhelv Universit\'e Pierre et Marie Curie -
Paris 6}
\vskip8pt
\centerline{\sevenhelv Tour 46 - 0\quad-\quad Bo\^\i te 186
\quad-\quad 4, place
Jussieu\quad-\quad 75252 PARIS
CEDEX 05}
\centerline{\sevenhelv T\'el : (33-1) 44 27 53 49
\quad-\quad T\'el\'ecopie : (33-1) 44 27 25 55\quad-\quad
lana@ccr.jussieu.fr}}
\null\kern-16mm
\line{\kern-10mm\raise-7pt\box10
\hfil\box12\hskip1cm}
\vskip6pt
\line{\kern-9mm\hrulefill}
\vskip3mm
}

\baselineskip=18pt
\def\n{\noindent}
\hfuzz=0,5cm

\centerline {  Analyse Fonctionnelle}

\centerline {\bf Espaces $L_p$ non commutatifs \`a valeurs
vectorielles}

\centerline {\bf  et
applications compl\`etement $p$-sommantes.}

Note de Gilles Pisier

Pr\'esent\'ee par
\vskip12pt

\underbar {R\'esum\'e franiais}. Soit $E$ un espace d'op\'erateurs au
sens de la
th\'eorie d\'evelopp\'ee r\'ecemment par Blecher-Paulsen et
Effros-Ruan. On introduit
une notion d'espace $L^p$ non commutatif \`a valeurs dans $E$ pour $1
\leq p <
\infty$ et on d\'emontre qu'elle poss\`ede les propri\'et\'es
naturelles que l'on
attend, pour la dualit\'e et l'interpolation par exemple. Cela permet
de d\'efinir
une notion d'application compl\`etement $p$-sommante adapt\'ee \`a la
cat\'egorie des
espaces d'op\'erateurs. Cette notion g\'en\'eralise celle introduite
pr\'ec\'edemment par
Effros-Ruan pour $p=1$.
\vskip12pt
\underbar {English title}. Noncommutative vector valued $L_p$-spaces
and
completely $p$-summing maps.

\underbar {English Abstract}. Let $E$ be an operator space in the sense
of the
theory recently developed by Blecher-Paulsen and Effros-Ruan. We
introduce a
notion of $E$-valued non commutative $L_p$-space for $1 \leq p <
\infty$ and we
prove that the resulting operator space satisfies the natural
properties to be
expected with respect to e.g. duality and interpolation. This notion
leads to
the definition of a ``completely p-summing" map which is the operator
space
analogue of the $p$-absolutely summing maps in the sense of
Pietsch-Kwapie\'n.
These notions extend the particular case $p=1$ which was previously
studied by
Effros-Ruan.

\vfill\eject
\centerline {\bf English Abridged version} We will work in the category
of
operator spaces as developped in the papers [1,2] and [4]. By an
operator space
we mean a closed subspace of the space $B(H)$ of all bounded operators
on $H$ for
some Hilbert space $H$. In this category the morphisms (resp.
isomorphisms) are
the completely bounded maps (resp. complete isomorphisms). We refer to
[1,2,4]
for more informations and in particular for the notions of
completely bounded map (in short c.b.), of  dual space and
quotient space in this category. We refer to our recent work [9,10]
for the
definition and basic properties of the complex interpolation method and
ultraproducts in the category of operator spaces, as well as for the
definition
and basic properties of the operator Hilbert space $OH(I)$ associated
to any
index set $I$.

Let $H,K$ be Hilbert spaces. We will denote by $H \otimes_2 K$ the
Hilbertian
tensor product. Let $E \subset B(H),\ F \subset B(K)$ be two operator
spaces. We
will denote by $E \otimes_m F$ the minimal (or spatial) tensor product,
\ie
the completion of the linear product $E \otimes F$ for the norm induced
by $B(H
\otimes_2 K)$.

We will denote by $S_p(K)$ the Schatten class formed of all the compact
operators $T : K \to K$ such that $tr |T|^p < \infty$, equipped with
the norm
$\N{T}_p = (tr |T|^p)^{1/p}$. We denote by $S_\infty (K)$ the class of
all
compact operators on $K$ equipped with the norm induced by $B(K)$.

We will define the space $S_p[K;E]$ when $E \subset B(H)$ is an
operator space.
First we define $S_\infty[K;E]-= S_\infty(K) \otimes_m E$ and
$S_1[K;E]-=
S_1(K) \otimes_\wedge E$ where $\otimes_\wedge$ is the operator space
version
of the projective tensor norm as introduced in [1,4] and developped by
Effros-Ruan in [6, 7, 8]. By these known results,
we have a
contractive inclusion $S_1[K;E]-\to S_\infty[K;E]$ which allows to
consider
this pair as a compatible couple in the sense of interpolation theory.
Then
for $1 < p < \infty$ and $\theta = 1/p$ we define the operator space
$S_p[K;E]$
as $S_p[K;E] = (S_\infty[K,E],S_1[K;E])_\theta$, where the operator
space
structure is defined as in [9].

To simplify our notation, let us restrict ourselves to the case $K =
\ell_2$.
In that case we will write $S_p$ and $S_p[E]$ instead of $S_p(K)$ and
$S_p[K;E]$. Assume $1 \leq p_0,p_1 \leq \infty$ and $p^{-1} =
(1-\theta)
p_0^{-1} + \theta p_1^{-1}$. Then $S_p[E] =
(S_{p_0}[E],S_{p_1}[E])_\theta$
completely isometrically. Moreover, let $p'=p(p-1)^{-1}$, then
$S_p[E]^* =
S_{p'}[E^*]$ completely isometrically. The space $S_2$ is completely
isometric
to the space $OH(\nat \times  \nat)$ in the sense of our previous work
[9].
Moreover $S_2[OH(I)]$ is completely isometric to $OH(\nat \times \nat
\times
I)$.

Let $K,L$ be Hilbert spaces, then we have the following analogue of
Fubini's
theorem : $S_p[K;S_p[L;E]] = S_p[K \otimes_2 L;E]-= S_p[L;S_p[K;E]]$
completely
isometrically.

Let $E,F$ be operator spaces. We denote by $cb(E,F)$ the space of all
c.b. maps
from $E$ into $F$ equipped with the c.b. norm.
Let $u : E \to F$ be a linear map. We say that $u$ is
completely $p$-summing if the operator $I_{S_p} \otimes u$ extends to a
bounded
map $\tilde u$ from $S_p \otimes_m E$ to $S_p[F]$. We denote
$\pi_p^0(u) = \N
{\tilde u}$.
Actually, it can be shown that $\tilde u$ is then completely bounded
with c.b.
norm $\N {\tilde u}_{cb} = \N {\tilde u}$. In particular $\N {u}_{cb}
\leq
\pi_p^0(u)$.

 We prove the analogue of the Pietsch factorization
theorem for these maps. Moreover we show that if $E \subset B(H)$ is
$n$
dimensional then $\pi_2^0(I_E) = n^{1/2}$. This implies that there is
an
isomorphism $u : E \to OH_n$ such that $\N {u}_{cb} \N {u^{-1}}_{cb}
\leq
n^{1/2}$ and a projection $P : B(H) \to E$ such that $\N {P}_{cb} \leq
n^{1/2}$.

For every operator $u : E \to OH(J)~~--(J$ an arbitrary set) we have
$\pi_2^0(u) = \pi_{2,oh}(u)$ where $\pi_{2,oh}(u)$ is the
$(2,oh)$-summing  norm
introduced in our preceding note [10]. Moreover, for every
 $u : OH(I) \to OH(J)~~-(I,J$ arbitrary sets) the  Hilbert Schmidt norm
  $\N {u}_{HS}$ of $u$ coincides with $\pi_2^0(u)$ and
$\pi_{2,oh}(u)$.
One can then reformulate a result of [10] as follows: A linear map $u :
E \to F$
is in $\Gamma_{oh}(E,F)$ iff there exists a constant $C$ such that for
all
 $n$ and all  $v : F \to OH_n$ we have
$\pi_2^0((vu)^*) \leq  C \pi_2^0(v)$.
 Equivalently, this means that for all $n$ and all
 $T :S_2^n \to S_2^n$ the norm
of $T \otimes u$ in $cb(S_2^n[E],S_2^n[F])$ is majorized by $C\N {T}$.
 Moreover, $\gamma_{oh}(u)$ is
equal to the smallest constant $C$ in either of
these properties. In particular this result yields a natural
operator space structure on the space
$\Gamma_{oh}(E,F)$, obtained by embedding it in a suitable direct sum
  (in the sense of $\ell^\infty)$ of copies of $cb(S_2[E],S_2[F])$.
\vfill\eject

Soit $1 \leq p < \infty$. On notera $S_p(K)$ l'espace des operateurs
compacts
$T : K \to K$ tels que $tr|T|^p < \infty$ muni de la norme $\N {T}_p =
(tr|T|^p)^{1/p}$. Pour $p=\infty$, on notera $S_\infty(K)$ l'ensemble
des
op\'erateurs compacts muni de la norme induite par $B(K)$. Si $K =
\ell_2$, on
note $S_p = S_p(\ell_2)$ et si $K = \ell_2^n$ on note $S_p^n =
S_p(\ell_2^n)$.

Il est bien connu que les espaces $S_p(K)$ (appel\'ees souvent classes
de
Schatten) sont un analogue non commutatif des espaces
$L_p(\Omega,\mu)$, tout
au moins pour un espace mesur\'e discret (\ie atomique). Dans le cas
commutatif,
pour tout espace de Banach $E$, on sait construire (suivant une id\'ee
attribu\'ee
\`a Bochner) l'espace $L_p(\Omega,\mu;E)$ des fonctions $L_p$ \`a
valeurs dans $E$.
Nous allons d\'efinir un analogue non commutatif $S_p[K;E]$ dans le cas
o\`u $E$ est
un espace d'op\'erateurs. L'espace $S_p[K;E]$ sera lui aussi un espace
d'op\'erateurs.

Pour simplifier nous ne consid\'erons dans cette note que le cas
discret.
N\'eanmoins, les id\'ees peuvent \^etre facilement adapt\'ees au cas
d'une alg\`ebre de
von Neumann $M$ munie d'une trace semifinie, fid\`ele et normale, \`a
condition que
l'alg\`ebre de von Neumann $M$ soit injective. Bien que les
d\'efinitions aient un
sens, des propri\'et\'es essentielles sont en d\'efaut si $M$ n'est pas
suppos\'ee
injective. Nous donnerons plus de d\'etails et les d\'emonstrations des
r\'esultats
annonc\'es ci-dessous dans une prochaine publication.

Soient $H,K$ deux espaces de Hilbert. On note $B(H)$ l'espace des
op\'erateurs
born\'es sur $H$ muni de sa norme usuelle. On appelle ``espace
d'op\'erateurs" un
sous-espace ferm\'e de $B(H)$. Dans la cat\'egorie des espaces
d'op\'erateurs telle
qu'elle est d\'evelop\'ee par Blecher-Paulsen [1,2] et Effros-Ruan [4],
les
morphismes (resp. isomorphismes, resp. isom\'etries) sont les
applications
compl\`etement born\'es (resp. les isomorphismes complets, resp. les
isom\'etries
compl\`etes). On abr\'egera compl\'etement born\'e en $c.b.$. On notera
$cb(E,F)$
l'espace des applications lin\'eaires de $E$ dans  $F$ muni de la norme
$\N
{~}_{cb}$. Nous renvoyons aux articles [1,2,4] pour les notions de dual
$E^*$
d'un espace d'op\'erateur $E$ et aussi pour la notion de quotient dans
cette
cat\'egorie. Nous renvoyons
\`a nos travaux pr\'ec\'edents [9,10] pour la d\'efinition   de la
m\'ethode d'interpolation complexe et des ultraproduits dans la
cat\'egorie des
espaces d'op\'erateurs, ainsi que pour la d\'efinition et les
principales propri\'et\'es
de l'espace d'op\'erateurs $OH(I)$ qui est l'analogue dans cette
cat\'egorie de
l'espace $\ell_2(I)$ associ\'e \`a un ensemble arbitraire $I$.

Soit $E \subset B(H),\ F \subset B(K)$ deux espaces d'op\'erateurs,
nous noterons
$E \otimes_m F$ leur produit tensoriel minimal (= spatial), \ie le
compl\'et\'e du
produit tensoriel alg\'ebrique $E \otimes F$ pour la norme induite par
$B(H
\otimes_2 K)$, o\`u on a not\'e $H \otimes_2 K$ le produit tensoriel
hilbertien de
$H$ et $K$. Nous renvoyons \`a [1, 2, 3, 5] pour la d\'efinition et les
principales
propri\'et\'es du produit tensoriel de Haagerup $E \otimes_h F$ de deux
espaces
d'op\'erateurs. Pour tout espace de Hilbert $K$, nous noterons $K_c$ et
$K_r$ les
espaces d'op\'erateurs (isom\'etriques \`a $K$) d\'efinis par $K_c =
B(\comp,K),K_r =
B(K^*,\comp)$.

Nous allons d\'efinir l'espace d'op\'erateurs $S_p[K;E]$. Commenions
par les cas $p
= \infty$ et $p = 1$. Pour $p = \infty$, on pose $S_\infty[K;E] =
S_\infty(K)
\otimes_m E$. Pour $p=1$, on d\'efinit $S_1[K;E] = S_1(K)
\otimes_\wedge E$  o\`u
$\otimes_\wedge$ d\'esigne l'analogue du produit tensoriel projectif
dans la
cat\'egorie des espaces d'op\'erateurs tel qu'il est d\'efini dans [1,
4]. Dans une
s\'erie d'articles Effros et Ruan [6, 7, 8] ont d\'evelopp\'e les
principales
propri\'et\'es de ce produit tensoriel. En particulier, d'apr\`es leurs
travaux il y
a une inclusion  contractive $S_1[K;E]-\to S_\infty[K;E]$.
Cette inclusion nous permet de consid\'erer le couple
$(S_1[K;E],S_\infty[K;E])$ comme un couple compatible pour
l'interpolation. On d\'efinit alors pour $1 < p < \infty$
l'espace $S_p[K;E]$ par la m\'ethode d'interpolation complexe
$$S_p[K;E]-= (S_\infty[K;E],S_1[K;E])_\theta$$ o\`u $\theta =
1/p$. La structure d'espace d'op\'erateurs sur $S_p[K;E]$ est
d\'efinie (voir [9]) par l'identit\'e  $$M_n(S_p[K;E]) =
(M_n(S_\infty[K;E]),M_n(S_1[K;E]))_\theta.$$ Rappelons que
l'on peut consid\'erer $(K_r,K_c)$ comme un couple compatible
pour l'interpolation grece \`a l'isom\'etrie lin\'eaire $i : K_c
\to K_r$ qui associe \`a tout $T : \comp \to K$ l'op\'erateur
adjoint $i(T) : K^* \to \comp$. On peut donc d\'efinir (voir
[9]) l'espace d'op\'erateurs $K_r(\theta) = (K_r,K_c)_\theta$
pour $0 < \theta < 1$. Nous poserons par convention $K_r(0)
= K_r$ et $K_r(1) = K_c$. Rappelons que l'on a les
isomorphismes compl\'etement isom\'etriques suivants (\cf [1,
5]). $$K_c \otimes_h E \otimes_h K_r = S_\infty(K)
\otimes_m E-\quad {\rm et}\quad-K_r \otimes_h E \otimes_h
K_c = S_1(K) \otimes_\wedge E.$$ D'apr\`es [9] on peut en
d\'eduire

\proclaim Th\'eor\`eme 1. Pour $1 < p < \infty$ et $\theta = 1/p$. On a
un
isomorphisme compl\`etement isom\'etrique
$$S_p[K;E] = K_r(1-\theta) \otimes_h E \otimes_h K_r(\theta).$$

\n En particulier, pour $\theta=1/2$ et (pour all\'eger la notation) si
$K=\ell_2$,
on a $S_2[E]= OH\otimes_h E \otimes_h OH$, compl\`etement
isom\'etriquement.

\n Voici
quelques propri\'et\'es de l'espace $S_p[K;E]$:

\n Soit $E,F$ des espaces
d'op\'erateurs et soit $u : E \to F$ une application $c.b.$. Alors
$I_{{S_p}(K)}
\otimes u$ s'\'etend en une application $U : S_p[K;E] \to S_p[K;F]$
telle que $\N
{U}_{cb} = \N {u}_{cb}$. Si $u$ est une isom\'etrie compl\`ete, il en
est de m\^eme
de $U$. Le produit tensoriel alg\'ebrique $S_p(K) \otimes E$ est dense
dans
$S_p[K;E]$. On a une formule de dualit\'e
$$S_p[K;E]^* = S_{p'}[K;E^*],~~~~~~~~~~({\rm o\grave u}~-p' =
p(p-1)^{-1})$$
compl\`etement isom\'etriquement. Quant \`a l'interpolation, soit
$(E_0,E_1)$ un
couple compatible d'espaces d'op\'erateurs, soit $E_\theta =
(E_0,E_1)_\theta$
l'espace d'op\'erateurs d\'ecrit dans [9], soit $1 \leq p_0,p_1,p \leq
\infty$ tels
que $p^{-1} = (1-\theta) p_0^{-1} + \theta p_1^{-1}$. On a alors
$$S_p[K,E_\theta] = (S_{p_0}[K;E_0],S_{p_1}[K;E_1])_\theta$$
compl\`etement isom\'etriquement.

On peut donner aussi l'analogue suivant du th\'eor\`eme de Fubini. Soit
$1 \leq p
\leq \infty$. Soit $K,L$ deux espaces de Hilbert. On a compl\`etement
isom\'etriquement
$$S_p[K;S_p[L;E]]-= S_p[K \otimes_2 L;E]-= S_p[L;S_p[K;E]].$$
Plus g\'en\'eralement, si $p \leq q \leq \infty$ on a une inclusion
compl\'etement
contractive $S_p[K;S_q[L;E]]$ dans $S_q[L;S_p[K;E]]$.
Soit $I$ le cardinal d'une base orthonormale de $K$. Alors $S_2(K)$ en
tant
qu'espace d'op\'erateurs est compl\`etement isom\'etriquement
identifiable \`a l'espace
$OH(I \times  I)$ introduit dans [9]. De plus pour tout ensemble $J$ on
peut
v\'erifier que $S_2[K;OH(J)]$ s'identifie compl\`etement
isom\'etriquement \`a $S_2(K)
\otimes_h OH(J)$ ou encore \`a $OH(I \times  I \times  J)$.

Le th\'eor\`eme suivant permet de ``calculer" la norme dans $S_p[K;E]$
et
$M_n(S_p[K;E])$.

\proclaim THEOREME 2. Soit $1 \leq p < \infty$. Soit $u \in S_p[K;E]$.
On a
$$\N {u}_{{S_p}[K;E]} = \inf \{\N {a}_{{S_{2p}(K)}} \N
{v}_{{S_\infty}[K;E]}\N
{b}_{{S_{2p}(K)}}\} \leqno (1)$$
o\`u l'infimum porte sur toutes les repr\'esentations de $u$ de la
forme
$u = (a \otimes I_E) v(b \otimes I_E)$, avec $a,b \in S_{2p}(K)$ et $v
\in
S_\infty[K;E]$.\sn
\n D'autre part, soit $F \subset B(L)$ un autre espace d'op\'erateurs.
Soit $a,b$ dans $S_{2p}[L]$. Notons $M(a,b)$ l'application de $B(L)$
dans
$S_p(L)$ d\'efinie par $M(a,b)y  = a y b$. On notera $\tilde M(a,b)$
l'application de $S_p[K;E]-\otimes_m B(L)$ dans $S_p[K;E]-\otimes_m
S_p(L)$
associ\'ee \`a $I \otimes M(a,b)$. Alors, pour tout $x$ dans $S_p[K,E]
\otimes_m
F$, sa norme $\N {x}_m$ est donn\'ee par
$$\N {x}_m = \sup \{\N {\tilde M(a,b) x}_{{S_p}[K \otimes_2 L;E]}\}
\leqno (2)$$
o\`u le supremum porte sur tous les $a,b$ dans la boule unit\'e de
$S_{2p}(L)$.

La formule (1) d\'ecrit $S_p[K;E]$ comme espace de Banach et (2)
d\'ecrit sa
structure comme espace d'op\'erateurs.

A l'aide des r\'esultats pr\'ec\'edents, on peut d\'evelopper une
th\'eorie des
applications $p$-sommantes entre espaces d'op\'erateurs tout-\`a-fait
analogue \`a
celle de Pietsch et Kwapie\'n (\cf e.g. [11]) pour les espaces de
Banach.

Soit $E,F$ deux espaces d'op\'erateurs.
Soit $u : E \to F$ une application lin\'eaire.
Nous dirons que $u$ est compl\`etement $p$-sommante si l'application
$\tilde u =
I_{S_p} \otimes u$ est born\'ee de $S_p \otimes_m E$ dans $S_p[F]$. On
pose
$\pi_p^0(u) = \N {\tilde u}_{{S_p} \otimes_m E \to S_p(F)}$. On note
$\Pi_p^0(E,F)$ l'espace des applications compl\`etement $p$-sommantes
et on le
munit de la norme $\pi_p^0$ pour laquelle c'est un espace de Banach.
Soit $v :
E_1 \to E$ et $w : F \to F_1$ des applications $c.b.$ entre espaces
d'op\'erateurs, on a alors
$$\pi_p^0 (w u v) \leq \N {w}_{cb} \pi_p^0(u) \N {v}_{cb}$$
De plus on peut noter que si $\tilde u$ est born\'e (\ie si $u$ est
compl\`etement
$p$-sommante) alors n\'ecessairement $\tilde u$ est $c.b.$ et $\N
{\tilde u}_{cb}
= \N {\tilde u} = \pi_p^0(u)$. En particulier on a
$$\N {u}_{cb} \leq \pi_p^0(u). \leqno (3)$$
Cette remarque permet de munir $\PI_p^0(E,F)$ de la structure d'espace
d'op\'erateurs naturellement induite par l'espace $cb (S_p \otimes_m
E,S_p[F])$.

On peut montrer par exemple que l'op\'erateur $M : B(\ell_2) \to S_p$
d\'efini par
$M(x) = a x b$ avec $a,b \in S_{2p}$ est un op\'erateur compl\`etement
$p$-sommant.
De plus, soit $E \subset B(\ell_2)$ et soit $E_p$ la fermeture de
$M(E)$ dans $S_p$. Alors l'application $M_1 : E \to E_p$ qui est la
restriction
de $M$ est compl\`etement
$p$-sommante.
Le th\'eor\`eme qui suit est l'analogue du th\'eor\`eme de
factorisation de Pietsch. Il
montre que l'exemple pr\'ec\'edent est fondamental.

\proclaim THEOREME 3. Soit $E \subset B(H)$. Soit $\tilde H = H \oplus
H
\oplus...$ la somme hilbertienne d'une famille d\'enombrable de copies
de $H$ et
soit $\pi : B(H) \to B(\tilde H)$ la repr\'esentation somme directe
d'une famille
d\'enombrable de copies de la repr\'esentation identique de $B(H)$.
Soit $u : E \to
F$ une application compl\`etement $p$-sommante et soit $C =
\pi_p^0(u)$. Il
existe alors un ensemble $I$ muni d'un ultrafiltre ${\cal U}$ et des
familles
$(a_\alpha)_{\alpha \in I},(b_\alpha)_{\alpha \in I}$ dans la boule
unit\'e de
$S_{2p}(\tilde H)$ telles que l'on ait pour tout $n$ et tout $(x_{ij})$
dans
$M_n(E)$
$$\N {(u(x_{ij}))}_{M_n(F)} \leq C \lim_{\cal U} \N {(a_\alpha
\pi(x_{ij})
b_\alpha)}_{M_n(S_p(\tilde H))}.\leqno (4)$$
R\'eciproquement, toute application v\'erifiant (4) pour tout $n$ est
n\'ecessairement compl\`etement $p$-sommante et telle que $\pi_p^0(u)
\leq C$.

\underbar {Remarques} : (i) Si $H$ est fini dimensionnel et si $p = 2$,
on peut
d\'emontrer le r\'esultat pr\'ec\'edent avec pour $\tilde H$ une somme
hilbertienne
d'un nombre fini $m\  (m \leq n^4 + 1)$ de copies de $H$. Dans ce cas
la boule
unit\'e de $S_4(\tilde H)$ \'etant compacte, on trouve (4) avec
$a_\alpha =
a~~-b_\alpha = b~~~-\forall~-\alpha \in I$ et avec $a,b$ dans la boule
unit\'e de
$S_4(\tilde H)$.

(ii) Dans le cas $p=2$, on trouve ainsi que si $u$ est compl\`etement
$2$-sommant, alors $u \in \Gamma_{oh} (E,F)$ au sens de [10] et
$\gamma_{oh}(u)
\leq \pi_2^0(u)$.

\n De plus si $E \subset B(H)$ alors $u : E \to F$ admet une extension
$\hat u :
B(H) \to F$ telle que $\pi_2^0(\hat u) = \pi_2^0(u)$. On peut montrer
que pour
tout sous espace $E \subset B(H)$ avec $\dim E = n$ on a $\pi_2^0(I_E)
=
n^{1/2}$. Il existe donc un isomorphisme $u : E \to OH_n$ tel que $\N
{u}_{cb}
\N {u^{-1}}_{cb} \leq \sqrt n$ et une projection $P : B(H) \to E$ telle
que $\N
{P}_{cb} \leq \sqrt n$.

(iii) Pour tout op\'erateur $u : E \to OH(J)~~--(J$ un ensemble
arbitraire) on a
$\pi_2^0(u) = \pi_{2,oh}(u)$ o\`u $\pi_{2,oh}(u)$ est la norme
d\'efinie dans notre
note pr\'ec\'edente [10].

(iv) Pour tout $u : OH(I) \to OH(J)~~-(I,J$ des ensembles arbitraires)
la norme
de Hilbert Schmidt $\N {u}_{HS}$ de $u$ coincide avec $\pi_2^0(u)$ et
$\pi_{2,oh}(u)$.

(v) On peut alors reformuler un r\'esultat de [10] de la mani\`ere
suivante : Une
application $u : E \to F$ est dans $\Gamma_{oh}(E,F)$ ssi il existe une
constante $C$ telle que pour tout $n$ et tout $v : F \to OH_n$ on a
$\pi_2^0((vu)^*) \leq  C \pi_2^0(v)$. Cela revient \`a dire que pour
tout $n$ et
tout $T :S_2^n \to S_2^n$ la norme de $T \otimes u$ dans
$cb(S_2^n[E],S_2^n[F])$ est major\'ee par $C\N {T}$. De plus
$\gamma_{oh}(u)$ est
\'egal \`a la plus petite constante $C$ dans l'une ou l'autre de ces
propri\'et\'es.
En particulier ce r\'esultat permet de munir l'espace
$\Gamma_{oh}(E,F)$ d'une
structure naturelle d'espace d'op\'erateurs en le plongeant dans une
somme
directe convenable (au sens $\ell^\infty)$ de copies de
$cb(S_2[E],S_2[F])$.
\vfill\eject

\centerline{\bf R\'ef\'erences bibligraphiques:}
\vskip12pt

 \item{[1]} D. Blecher and V. Paulsen. Tensor products of operator
 spaces.
 J. Funct. Anal. 99 (1991) 262-292.

\item{[2]} D. Blecher. Tensor products of operator spaces II.
(Preprint)
 1990.   Canadian J. Math. A para\^{\i}tre.

 \item{[3]} D. Blecher and R. Smith. The dual of the Haagerup tensor
 product. Journal London Math. Soc. 45 (1992) 126-144.

 \item{[4]} E. Effros and Z.J. Ruan. A new approach to operators
 spaces.
 Canadian Math. Bull.
34 (1991) 329-337.

 \item{[5]}  E. Effros and Z.J. Ruan. Self duality for the Haagerup
 tensor product and Hilbert space factorization.  J. Funct. Anal. 100
(1991) 257-284.

 \item{[6]}  E. Effros and Z.J. Ruan. Mapping spaces
and liftings for operator spaces. (Preprint) Proc. London
Math. Soc.  A para\^{\i}tre.

\item{[7]}  E. Effros and Z.J. Ruan. The
Grothendieck-Pietsch and Dvoretzky-Rogers Theorems for
operator spaces. (Preprint 1991) J. Funct. Anal. A para\^{\i}tre.

\item{[8]}  E. Effros and Z.J. Ruan. On approximation
properties for operator spaces, International J. Math. 1
(1990) 163-187.

\item {[9]} G.  Pisier.  Espace de Hilbert
d'op\'erateurs et interpolation complexe. Comptes Rendus
Acad. Sci. Paris S\'erie I, 316 (1993)  47-52.

\item {[10]} G.  Pisier. Sur les
op\'erateurs factorisables par $OH$.  Comptes Rendus Acad.
Sci. Paris S\'erie I, 316 (1993) 165-170.

\item {[11]} G.  Pisier. Factorization of linear
operators and the Geometry of Banach spaces.  CBMS
(Regional conferences of the A.M.S.)    60, (1986),
Reprinted with corrections 1987.

\vskip24pt
Universit\'e Paris 6

Equipe d'Analyse, Bo\^\i te 186,

75252 Paris Cedex 05, France

 \end

\end

For every operator $u : E \to OH(J)~~--(J$ an arbitrary set) we have
$\pi_2^0(u) = \pi_{2,oh}(u)$ where $\pi_{2,oh}(u)$ is the
$(2,oh)$-summing  norm
introduced in our preceding note [9]. Moreover, for every
 $u : OH(I) \to OH(J)~~-(I,J$ arbitrary sets) the  Hilbert Schmidt norm
  $\N {u}_{HS}$ of $u$ coincides with $\pi_2^0(u)$ and
$\pi_{2,oh}(u)$.
One can then reformulate a result of [9] as follows: A linear map $u :
E \to F$
is in $\Gamma_{oh}(E,F)$ iff there exists a constant $C$ such that for
all
 $n$ and all  $v : F \to OH_n$ we have
$\pi_2^0((vu)^*) \leq  C \pi_2^0(v)$.
 Equivalently, this means that for all $n$ and all
 $T :S_2^n \to S_2^n$ the norm
of $T \otimes u$ in $cb(S_2^n[E],S_2^n[F])$ is majorized by $C\N {T}$.
 Moreover, $\gamma_{oh}(u)$ is
equal to the smallest constant $C$ in either of
these properties. In particular this result yields a natural
operator space structure on the space
$\Gamma_{oh}(E,F)$, obtained by embedding it in a suitable direct sum
  (in the sense of $\ell^\infty)$ of copies of $cb(S_2[E],S_2[F])$.